\documentclass[11pt]{amsart}
\usepackage{graphicx}
\allowdisplaybreaks
\usepackage{amsmath,amssymb,mathrsfs,amsthm}

\usepackage{verbatim}
\usepackage[spanish,USenglish]{babel} % espanol, ingles
\usepackage{color}
\usepackage{tikz}
\usepackage{graphicx}

\newtheorem{thm}{Theorem}[section]

\newtheorem{proposition}[thm]{Proposition}
\newtheorem{theorem}[thm]{Theorem}
\newtheorem{lemma}[thm]{Lemma}

\theoremstyle{definition}
\newtheorem{definition}[thm]{Definition}

\theoremstyle{remark}
\newtheorem{remark}[thm]{Remark}

\numberwithin{equation}{section}

\def\N{\mathbb{N}}

\def\C{\mathbb{C}}

\def\Z{\mathbb{Z}}

\def\P{\mathcal{P}}

\def\Tn2{{\mathcal T}^{(n)}_2}

\def\N{\mathbb{N}}
\def\C{\mathbb{C}}

\begin{document}

\title[Powers of Catalan generating functions ]
{Powers of Catalan generating functions for bounded operators}

\author[Miana]{Pedro J. Miana}
\address{Departamento de Matem\'aticas, Instituto Universitario de Matem\'aticas y Aplicaciones, Universidad de Zaragoza, 50009 Zaragoza, Spain.}
\email{pjmiana@unizar.es}

\author[Romero]{Natalia Romero}
\address{Departamento de Matem\'aticas y Computaci\'on,  Universidad de la Rioja, 26006 Logro\~{n}o, Spain.}
\email{natalia.romero@unirioja.es}

\thanks{  Pedro J. Miana has been partially supported by Project ID2019-105979GBI00, DGI-FEDER, of the MCEI and Project E48-20R, Gobierno de Arag\'on, Spain. Natalia Romero   has been partially supported by the project MTM2018-095896-B-C21
of the Spanish Ministry of Science. \newline Pedro J. Miana, corresponding author, Departamento de Matem\'aticas, Instituto Universitario de Matem\'aticas y Aplicaciones, Universidad de Zaragoza, 50009 Zaragoza, Spain.pjmiana@unizar.es. \newline Natalia Romero, Departamento de Matem\'aticas y Computaci\'on,  Universidad de la Rioja, 26006 Logro\~{n}o, Spain. natalia.romero@unirioja.es}



\maketitle

\begin{abstract} Let $c=(C_n)_{n\ge 0}$ be the Catalan sequence and $T$ a linear and bounded operator on a Banach space $X$ such $4T$  is a power-bounded operator. The Catalan generating function is defined by the following Taylor series,
$$
C(T):=\sum_{n=0}^\infty C_nT^n.
$$
Note that the operator $C(T)$ is a solution of the quadratic equation $TY^2-Y+I=0.$  In this paper we define  powers of the Catalan generating function $C(T)$ in terms of the Catalan triangle numbers. We obtain new formulae which involve Catalan triangle numbers; the spectrum of $c^{\ast j}$  and the expression of  $c^{-\ast j}$ for $j\ge 1$ in terms of Catalan polynomials ($\ast$ is the usual convolution product in sequences). In the last section, we give some particular examples to illustrate our results and some ideas to continue this research in the future. \end{abstract}

%%%%%%%%%%%%%
\section{Introduction}
%%%%%%%%%%%%%
The well-known Catalan numbers $(C_n)_{n\ge 0}$ are  given by the combinatorial formula
$$
 C_n={1\over n+1}{2n\choose n},\quad  \ n\ge 0,
$$
They may be defined recursively by $C_0=1$ and
\begin{equation}\label{cata}C_n= \sum_{i=0}^{n-1} C_i C_{n-1-i}, \qquad n\ge 1, \end{equation}
 and first terms in this sequence are
$
1,\,\,1,\,\,2,\,\,5,\,\,14, \,\,42, \,\,132, \dots .
$
They appear in a wide range of combinatorial problems: they count the  number of ways to triangulate a regular polygon with $n +2$ sides; or,
the number of  ways  that   $2n$ people seat around a circular table are simultaneously
shaking hands with another person at the table in such a way that none of the arms cross each other, see for example \cite{[Sl1], [St2]}.

The generating function of the Catalan sequence $c=(C_n)_{n\ge 0}$ is defined by
\begin{equation}\label{gene}
C(z):=\sum_{n=0}^\infty C_nz^n= {1- \sqrt{1-4z}\over 2z}, \quad z\in D(0,{1\over 4}):=\{z\in \C\,\,|\,\, \vert z\vert<{1\over 4}\}.
\end{equation}
This function satisfies the quadratic equation
$
zy^2-y+1=0
$.

The main aim in \cite{[MR3]} is to consider the quadratic equation
\begin{equation}\label{Ceq}
TY^2-Y+I=0,
\end{equation}
 in the set of linear and bounded operators, ${\mathcal B}(X)$ on a Banach space $X$, where $I$ is the identity on the Banach space, and  $T , Y \in {\mathcal B}(X)$.  Formally, some solutions of this vector-valued quadratic equations are expressed by
$$
Y={1\pm\sqrt{1-4T}\over 2T},
$$
 which involves the (non-trivial) problems of the square root of operator $1-4T$ and the inverse of operator $T$.

In this paper, we are concerned about the powers of $(C(T))^n$ for $n\in \Z$ and it is organized as follows. In the second section we consider the Catalan triangle sequences $(B_{n,k})_{n\ge 1, 1\le k \le n}$ and $(A_{n,k})_{n\ge 1, 1\le k \le n+1}$. We prove new formulae for these numbers (Lemma \ref{tec}) and their asymptotic estimation (Lemma \ref{assym}). We treat polynomials and generating formulae for these Catalan triangle numbers, see Definition \ref{polys} and Theorem \ref{generating}.

In third section, we consider the Banach algebra $(\ell^1(\N^0, {1\over 4^n}), \Vert \quad\Vert_{1, {1\over 4^n}}, \ast)$, where $$
\Vert a\Vert_{1, {1\over 4^n}}:=\sum_{n=0}^\infty {\vert a(n)\vert\over 4^n}<\infty, \qquad
(a\ast b)(n)=\sum_{j=0}^na{(n-j)}b(j), \quad n\ge 0,
$$
where  $a,b \in \ell^1(\N^0, {1\over 4^n}).$ We consider Catalan triangle sequences $(a_k)_{k\ge 1},$ $(b_k)_{k\ge 1}$ $\subset \ell^1(\N^0, {1\over 4^n})$ (Definition \ref{suc}). These sequences are powers of the Catalan sequence $c$ in $\ell^1(\N^0, {1\over 4^n})$ (Proposition  \ref{cons}); we describe their spectrum set in Proposition \ref{spec}. An original and motivating results connects $c^{-\ast k}$ and Catalan polynomials in Theorem \ref{convoinve}.

The powers of the Catalan generating operator $C(T)$ are studied in fourth section. We transfer our results from the algebra $\ell^1(\N^0, {1\over 4^n})$ to ${\mathcal B}(X)$ via  the algebra homomorphism $\Phi$,
$$
\Phi(a)x:=\sum_{n\ge 0} a_nT^n(x), \qquad  a=(a_n)_{n\ge 0}\in  \ell^1(\N^0, {1\over 4^n}) , \quad x\in X,
$$
 Note that $\Phi(c)=C(T)$, $\Phi(b_k)=(C(T))^{2k}$ and $\Phi(a_k)=(C(T))^{2k-1}$ for $k\ge 1$. We  describe $(C(T))^{-j}$ in terms of Catalan polynomials; we estimate their norms and describe $\sigma ((C(T))^{j})$ for $j\in \Z$ in Theorem \ref{genne}.

In the last section we illustrate our results with some  concrete operators $T$ in  the equation  (\ref{Ceq}). We consider the Euclidean space $\C^2$ and matrices $$T= \begin{pmatrix}
\lambda & 0\\
0 & \mu
\end{pmatrix}, \,\,\begin{pmatrix}
0 & \lambda\\
\lambda & 0
\end{pmatrix},\,\,
 \begin{pmatrix}
\lambda & \mu\\
0 & \lambda
\end{pmatrix}.$$
We solve the equation (\ref{Ceq}) and  calculate $(C(T))^{j}$ for these matrices and $j\in \Z$. We also check $(C(a))^{j}$ for some particular values of $a\in \ell^1(\N^0, {1\over 4^n})$ and $j\ge 1$. Finally we present some ideas to continue this research.

\section{ Some new results about Catalan triangle numbers}\label{GCHO}
%%%%%%%%%%%%%%%%%%%
%%%%%%%%%%%%%%%%%%%%

Calatan triangle numbers $(B_{n,k})_{n\ge 1, 1\le k\le n}$ were introduced in \cite{[Sh74]}. These combinatorial numbers $
B_{n,k}$
are the  entries of the following Catalan triangle:
\begin{equation} \label{shapiro}
\begin{tabular}{c|ccccccc}
  $n\setminus k$ &1  & 2 & 3 & 4 & 5&6&\dots
  \\
    \hline
  1 & 1 & \ &  \ &  \ &\ & \ & \\
  2 & 2 & 1& \ & \ &\  &\  & \ \\
  3 & 5 &  4& 1&\  &\  & \ & \ \\
  4 & 14 & 14 & 6& 1 &\ &\ & \ \\
  5 & 42 & 48 & 27 & 8&1 &\ & \ \\
  6 & 132 & 165 & 110 & 44&10 &1 & \  \\
  \dots & \dots & \dots & \dots & \dots&\dots &\dots & \dots \\
\end{tabular}
\end{equation}
which are given by
\begin{equation}\label{numbers1}B_{n,k}:={k\over n}{2n\choose n-k}, \ n,k\in \N,\ k\le n.
\end{equation}
Numbers $B_{n,k}$ has several applications: they count the number of leaves at level $k+1$ in all ordered trees with $n+1$ edges; $B_{n,k}$ is also the number of walks of $n$ steps, each in direction $N$, $S$, $W$ or $E$, starting at the origin, remaining
in the upper half-plane and ending at height $k$; or $B_{n, k}$ denotes the number of pairs of non-intersecting paths of length $n$ and distance $k$, see for example \cite{[Sh74]} and sequence A039598 in  \cite{[Sl]}. Notice that $B_{n,1}=C_n$ and $B_{n,n}=1$ for  $ n\ge 1$.

In the last years,  Catalan triangle (\ref{shapiro}) has been studied in detail.
These numbers $(B_{n,k})_{ n\ge  k \ge 1 }$ have been analyzed
in many ways. For instance,  symmetric functions have been used in \cite{[CC]},  recurrence relations  in \cite{[Sla]}, or in \cite{[GZ]}  the Newton interpolation formula, which is applied to conclude  divisibility properties of sums of products of binomial coefficients.

Other  combinatorial numbers $A_{n,k}$
defined as follows
\begin{equation}\label{numbers2}
A_{n,k}:={2k-1\over 2n+1}{2n+1\choose n+1-k}, \ n,k\in \N,\ k\le n+1,
\end{equation}
appear as the entries of this other Catalan triangle,
\begin{equation}\label{shapiro2}
\begin{tabular}{c|ccccccc}
  $n\setminus k$ &1  & 2 & 3 & 4 & 5&6&\dots
  \\
    \hline
   0 & 1 & \ &  \ &  \ &\ & \ & \\
   1 & 1 & 1 &  \ &  \ &\ & \ & \\
  2 & 2& 3 & 1& \  &\  &\  & \ \\
  3 & 5& 9 &  5& 1&\    & \ & \ \\
  4 & 14& 28 &20 & 7& 1 &\  & \ \\
  5 & 42 & 90 & 75& 35&9&1 & \ \\
  6 & 132 & 297 & 275 & 154&54 &11 & 1  \\
  \dots & \dots & \dots & \dots & \dots&\dots &\dots & \dots \\
\end{tabular}
\end{equation}
which  is considered in \cite{[MR2]}. These numbers also admit combinatorial interpretations: they count the number of lattice paths ending at a given height, in particular certain Grand Dyck paths, see more details in \cite{[G]} and sequence A039599 in  \cite{[Sl]}. Notice that $A_{n,1}=C_n$ and $A_{n,n+1}=1$ for $n\ge 1$.

The entries $
B_{n,k}$ and $A_{n,k}$ of the above two particular  Catalan triangles  satisfy  the  recurrence relations
\begin{equation}\label{recurr1}B_{n,k}=B_{n-1,k-1}+2B_{n-1,k}+B_{n-1,k+1},\qquad k\ge 2,\end{equation}
and
\begin{equation}\label{recurr2}A_{n,k}=A_{n-1,k-1}+2A_{n-1,k}+A_{n-1,k+1},\qquad k\ge 2.\end{equation}

The generating function of the Catalan sequence $(C_n)_{n\ge 0}$ is defined by
\begin{equation}\label{gene}
C(z):=\sum_{n=0}^\infty C_nz^n= {1- \sqrt{1-4z}\over 2z}, \quad z\in D(0,{1\over 4}):=\{z\in \C\,\,|\,\, \vert z\vert<{1\over 4}\}.
\end{equation}
Note that $C({1\over 4})=2$.

\begin{theorem}\label{resol} Take $z\in D(0,{1\over 4})$.
\begin{itemize}
\item[(i)] For $\lambda \not= C(z)$,
$$
{1 \over \lambda-C(z)}={\lambda z-1+zC(z)\over \lambda^2z-\lambda +1}.
$$

\item[(ii)] For $w\in D(0,{1\over 4})$ and $w\not={z\over (1+z)^2}$,
$$
{C^2(w)\over 1-zwC^2(w)}= {C(w)-(z+1)\over w(1+z)^2-z}.
$$
\end{itemize}
\end{theorem}

\begin{proof}(i) Note that
$$
(\lambda-C(z))(\lambda z-1+zC(z))=z\lambda^2-\lambda+C(z)-zC^2(z)=z\lambda^2-\lambda+1,
$$
for $\lambda \in \C$.

\item[(ii)]By item (i), we get that
\begin{eqnarray*}
{C^2(w)\over 1-zwC^2(w)}&=&{C^2(w)\over z}{1\over {1+z\over z}-C(w)}=C^2(w){w(1+z)-z+wzC(w)\over w(1+z)^2-z}\cr&=&{C(w)-1\over w}{w(1+z)-z+wzC(w)\over w(1+z)^2-z}= {C(w)-(z+1)\over w(1+z)^2-z},
\end{eqnarray*}
where we have applied again the equality $wC^2(w)-C(w)+1=0$. \end{proof}
As the following identity holds,
$$
C(z)^q=\sum_{n\ge 0}{q\over n+q}{2n-1+q\choose n}z^n, \quad q\ge 1, \quad z\in D(0,{1\over 4}),
$$
(\cite[Exercise A.32(a)]{[St2]}), we  take $q=2k $ and $q=2k+1$ for $k\ge 1$ to obtain the generating functions for the columns of the Catalan triangles,
\begin{eqnarray}
\sum_{n=k}^\infty B_{n,k}z^n&=&z^kC^{2k}(z)=(C(z)-1)^k,\qquad k\ge 1, \label{gene2} \\
\sum_{n=k}^\infty A_{n,k+1}z^n&=&z^{k}C^{2k+1}(z)=C(z)(C(z)-1)^k, \qquad k\ge 0, \label{gene3}
\end{eqnarray}
for  $ z\in D(0,{1\over 4})$.  Note that to get the second equality in both lines, we use the equality $zC^2(z)=C(z)-1$.

 We apply the formula (\ref{gene}) to get  $$\lim_{z\to {1\over 4}}C(z)=2,\qquad \lim_{z\to -{1\over 4}}C(z)=2(\sqrt{2}-1),$$(\cite[Exercise A.66]{[St2]}). Also other direct applications of Abel's theorem allows us to prove the following result.

\begin{lemma}\label{tec} Given $k\ge 1$,
\begin{eqnarray*}
\sum_{n=k}^\infty B_{n,k}{1\over 4^n}&= & 1,\qquad \sum_{n=k}^\infty B_{n,k}{(-1)^n\over 4^n}= (2\sqrt{2}-3)^k, \\
\sum_{n,k\ge 1}^\infty B_{n,k}{1\over 4^{n+k}}&= & {1\over 3},\qquad \sum_{n,k\ge 1}^\infty B_{n,k}{(-1)^n\over 4^{n+k}}= {8\sqrt{2}-13\over 41}, \\
\sum_{n=k}^\infty A_{n,k+1}{1\over 4^n}&= & 2, \qquad \sum_{n=k}^\infty A_{n,k+1}{(-1)^n\over 4^n}= 2(\sqrt{2}-1)(2\sqrt{2}-3)^k,\\
\sum_{n,k\ge 0}^\infty A_{n,k+1}{1\over 4^{n+k}}&= & {8\over 3},\qquad \sum_{n,k\ge 0}^\infty  A_{n,k+1}{(-1)^n\over 4^{n+k}}= {8\over 41}(5\sqrt{2}-3). \\
\end{eqnarray*}
\end{lemma}

\begin{proof} We apply formulae (\ref{gene2}) and (\ref{gene3}) in the points $z={1\over 4}$ and ${-1\over 4}$.
\end{proof}

In the next lemma, we extend the asymptotic estimation for Catalan numbers
$$
C_n\sim{4^n\over \sqrt{\pi}n^{3\over 2}}, \qquad n\to \infty,
$$
(\cite[Exercise A.64]{[St2]}) to Catalan triangle numbers.

\begin{lemma} \label{assym} Given $k\ge 1$,
\begin{eqnarray*}
B_{n,k}&\sim & {4^{n}\over \sqrt{\pi}}{k\over n^{3\over 2}}, \qquad n\to \infty,\\
A_{n,k}&\sim & {4^{n}\over \sqrt{\pi}}{2k-1\over n^{3\over 2}}, \qquad  n\to \infty.
\end{eqnarray*}
\end{lemma}
\begin{proof} We use the well-known Stirling formula $n!\sim e^{-n}n^n \sqrt{2\pi n}$ to show both equivalences.
\end{proof}

We now introduce the generating functions for the rows of the Catalan triangle numbers.

\begin{definition}\label{polys} Given $n\ge 0$, we define the polynomials
$$
P_n(z):=\sum_{j=0}^{n} B_{n+1, j+1}z^j,\qquad \qquad
Q_n(z):=\sum_{j=0}^{n+1} A_{n+1, j+1}z^j.
$$
\end{definition}

The first values of these families of polynomials are given by
\begin{eqnarray*}
P_0(z)=1,\qquad &\quad&  Q_0(z)=1+z \\
P_1(z)=2+z,\qquad &\quad&  Q_1(z)=2+3z+1 \\
P_2(z)=5+4z+z^2,\qquad &\quad& Q_2(z)=5+9z+5z^2+z^3 \\
P_3(z)=14+14z+6z^2+z^3,\qquad &\quad& Q_3(z)=14+28z+20z^2+7z^3+z^4 \\
\end{eqnarray*}

\begin{theorem}\label{recurre} (i)  The only  solution of the recurrence system
$$
 \left\{
 \begin{array}{lll}
R_0(z)=1,\qquad \\[1ex]
zR_n(z)+C_n  = (z+1)^2R_{n-1}(z),\quad n\ge 1,
  \end{array}
 \right.
$$
is the polynomial sequence $(P_n)_{n\ge 0}$ given in Definition \ref{polys}.

\medskip
(ii)  The only  solution of the recurrence system
$$
 \left\{
 \begin{array}{lll}
R_0(z)=1+z,\qquad \\[1ex]
zR_n(z)+C_n  = (z+1)^2R_{n-1}(z),\quad n\ge 1,
  \end{array}
 \right.
$$
is the polynomial sequence $(Q_n)_{n\ge 0}$ given in Definition \ref{polys}.

\end{theorem}

\begin{proof} It is enough to check that the sequence $(P_n)_{n\ge 0}$  satisfies the recurrence relation. Similarly the  polynomial sequence $(Q_n)_{n\ge 0}$ does. By the recurrence relation \ref{recurr1}, we get
\begin{eqnarray*}
P_{n+1}(z)&=&\sum_{j=0}^{n+1}B_{n+2,j+1}z^j=\sum_{j=0}^{n+1}(B_{n+1,j}+2B_{n+1,j+1}+B_{n+1,j+2})z^j\cr
&=&z\sum_{j=1}^{n+1}B_{n+1,j}z^{j-1}+2\sum_{j=0}^{n+1}B_{n+1,j+1}z^j+{1\over z}\sum_{j=0}^{n+1}B_{n+1,j+2}z^{j+1}\cr
&=&(z+2)P_n(z)+{1\over z}\left(\sum_{j=0}^{n}B_{n+1,j+1}z^{j}-B_{n+1,1}\right)\cr
&=& {(z+1)^2\over z}P_n(z)-{C_{n+1}\over z},
\end{eqnarray*}
and we conclude the equality.
\end{proof}

\begin{remark}\label{seq} {\rm The sequences of  polynomials  $(P_n)_{n\ge 0}$ and  $(Q_n)_{n\ge 0}$ are useful to prove equalities for Catalan triangles numbers and other sequences of integer numbers. For example,  taking $z=1$ in Theorem \ref{recurre}, we  prove easily by induction method that
$$
\sum_{k=1}^n B_{n,k}= {n+1\over 2}C_n, \qquad  \sum_{k=1}^{n+1} A_{n,k}= {(n+1)}C_n, \qquad n\ge 1.
$$
Indeed, we claim that $P_{n-1}(1)={n+1\over 2}C_n,$ for $n\ge 1$. As $P_n(1)+C_n  = 2^2P_{n-1}(1)=2(n+1)C_n$, we have that
$$
P_n(1)= (2n+1)C_n={n+2\over 2}C_{n+1},
$$
and we conclude the proof. An alternative proof appears in \cite[Proposition 3.1]{[Sh74]}. Similarly for $z=-1$, we get that
$$
\sum_{k=1}^n(-1)^k B_{n,k}= -C_{n-1}, \qquad  \sum_{k=1}^{n+1}(-1)^k A_{n,k}=0, \qquad n\ge 1.
$$
see, for example \cite[Theorem 2.1 and 2.2]{[MRInter]} and references therein.

 For $z={1\over 4}$, we follow similar ideas by induction method to get that
$$
\sum_{k=1}^n B_{n,k}\left({1\over 4}\right)^{k}={a(n)\over 4^n}, \qquad
\sum_{k=1}^{n+1} A_{n,k}\left({1\over 4}\right)^{k}={b(n)\over 4^{n+1}},
$$
where $(a(n))_{n\ge 1}$ is the integer sequence A194725 and $(b(n))_{n\ge 0}$ is A130970 given in the The On-Line Encyclopedia of Integer Sequences by N.J.A. Sloane, \cite{[Sl]}.

Finally for $z={-1\over 4}$, we obtain that
$$
\sum_{k=1}^n B_{n,k}\left({-1\over 4}\right)^{k}=-{d(n)\over (-4)^n}, \qquad
\sum_{k=1}^{n+1} A_{n,k}\left({-1\over 4}\right)^{k}=-{e(n)\over 4^{n+1}},
$$
where $(d(n))_{n\ge 1}$ is the integer sequence A051550 and $(e(n))_{n\ge 0}$ is A132863 given in \cite{[Sl]}.

}
\end{remark}

In the next theorem, we obtain the generating function for polynomial $(P_n)_{n\ge 0}$ and  $(Q_n)_{n\ge 0}$ given in Definition \ref{polys}.

\begin{theorem}\label{generating} For $n\ge 0$,
\begin{eqnarray*}
P(z,w):=\sum_{n\ge 0}P_n(z)w^n&=& {C(w)-(z+1)\over w(1+z)^2-z},\cr
Q(z,w):=\sum_{n\ge 0}Q_n(z)w^n&=& {(C(w)-(z+1))(z+1)\over w(1+z)^2-z}=P(z,w)(z+1).\cr
\end{eqnarray*}

\end{theorem}
\begin{proof} We take $z, w\in \C$ such that the bivariate generating function for polynomial $(P_n)_{n\ge 0}$ converges. Then
\begin{eqnarray*}
P(z,w)&=&\sum_{n\ge 0}P_n(z)w^n= \sum_{n\ge 0}\sum_{j=0}^{n} B_{n+1, j+1}z^jw^n=\sum_{j\ge 0}z^j\sum_{n=j}^{\infty} B_{n+1, j+1}w^n\cr
&=&\sum_{j\ge 0}z^jw^jC^{2j+2}(w)={C^2(w)\over 1-zwC^2(w)}= {C(w)-(z+1)\over w(1+z)^2-z},
\end{eqnarray*}
where we have applied the equation (\ref{gene2}), and Theorem \ref{resol} (ii).

Similarly,
\begin{eqnarray*}
Q(z,w)&=&\sum_{n\ge 0}Q_n(z)w^n= \sum_{n\ge -1}\sum_{j=0}^{n+1} A_{n+1, j+1}z^jw^n-{1\over w}\cr
&=&\sum_{j\ge 0}z^j\sum_{n=j-1}^{\infty} A_{n+1, j+1}w^n-{1\over w}
=\sum_{j\ge 0}z^jw^{j-1}C^{2j+1}(w)-{1\over w}\cr&=&{1\over w}{C(w)-1+zwC^2(w)\over 1-zwC^2(w)}=
{(1+z)C^2(w)\over 1-zwC^2(w)}\cr
&=&{(C(w)-(z+1))(z+1)\over w(1+z)^2-z}=P(z,w)(z+1),
\end{eqnarray*}
where we have applied the equation (\ref{gene3}), and   Theorem \ref{resol} (ii).
\end{proof}

\begin{remark} Note that for $\vert w\vert\le {1\over 4}$ and $\vert z\vert <1$, functions $P(z,w)$ and $Q(z, w)$ are well-defined due to
$$
\vert P(z,w)\vert \le \sum_{n\ge 0}\vert P_n(z)\vert {1\over 4^n}= 4\sum_{j\ge 0}\vert z\vert^j=4{1\over 1-\vert z\vert}.
$$

Formulae given in Theorem \ref{generating} extend several known generating formula, for example, for Catalan numbers
\begin{eqnarray*}
P(0,w) &= &\sum_{n\ge 0}P_n(0)w^n=\sum_{n\ge 0}B_{n+1,1}w^n= \sum_{n\ge 0}C_{n+1}w^n= {C(w)-1\over w},\cr
Q(0,w) &= &\sum_{n\ge 0}Q_n(0)w^n=\sum_{n\ge 0}A_{n+1,1}w^n= \sum_{n\ge 0}C_{n+1}w^n= {C(w)-1\over w}.
\end{eqnarray*}
Other generating functions for integer natural sequences, see  Remark \ref{seq}, are also obtained.

\end{remark}

\section{ Sequences of Catalan triangle numbers}\label{GCHO}

In this section, we consider the weighted Banach algebra $\ell^1(\N^0, {1\over 4^n})$.   This algebra is formed by sequences $a=(a(n))_{n\ge 0}$ such that
$$
\Vert a\Vert_{1, {1\over 4^n}}:=\sum_{n=0}^\infty {\vert a(n)\vert\over 4^n}<\infty,
$$
and the product is the usual convolution $\ast $ defined by
$$
(a\ast b)(n)=\sum_{j=0}^na{(n-j)}b(j), \qquad a,b \in \ell^1(\N^0, {1\over 4^n}).
$$
For $a, b\in \ell^1(\N^0, {1\over 4^n})$, note that
\begin{eqnarray*}
\Vert a\ast b\Vert_{1, {1\over 4^n}}&=&\sum_{n=0}^\infty {\vert (a\ast b)(n)\vert\over 4^n}\le \sum_{n=0}^\infty {1\over 4^n}\sum_{j=0}^n\vert a{(n-j)}\vert \,\vert b(j)\vert \cr &=&\sum_{j=0}^\infty \vert b(j)\vert \sum_{n=j}^\infty{\vert a{(n-j)}\vert \over 4^n}=\Vert a\Vert_{1, {1\over 4^n}}\, \Vert b\Vert_{1, {1\over 4^n}}.
\end{eqnarray*}
We write $a^{\ast 0}= a$ and $a^{\ast n}= a^{\ast(n-1)}\ast a$ for $n\in \N$.

The canonical base $\{\delta_j\}_{j\ge 0}$ is formed by sequences such that $(\delta_j)(n):=\delta_{j,n}$ is the known delta Kronecker. Note that $\delta_1^{\ast n}=\delta_n$ for $n\in \N$.
This Banach algebra has identity element, $\delta_0$, its spectrum set is the closed disc $\overline{D(0,{1\over 4})}$  and its Gelfand transform is given by the $Z$-transform
$$
Z(a)(z):= \sum_{n=0}^\infty {a(n)}z^n, \qquad z\in \overline{D(0,{1\over 4})},
$$
(\cite[Example 14.35]{Mu}). It is straightforward to check that $Z(\delta_n)(z)=z^n$ for $n\ge 0$ (see, for example, \cite[p. 21-22]{La}).

We recall that the  resolvent set of $a\in \ell^1(\N^0, {1\over 4^n})$, denoted as $\rho(a)$, is defined by $$\rho(a):=\{\lambda \in \C\, \, : \,\, (\lambda\delta_0-a)^{-1}\in \ell^1(\N^0, {1\over 4^n})\},
$$
and the spectrum set of  $a$  is denoted by $\sigma(a)$ and given by $\sigma(a):=\C\backslash \rho(a)$.

The Catalan numbers may be defined recursively by $C_0=1$ and
\begin{equation}\label{cata}C_n= \sum_{i=0}^{n-1} C_i C_{n-1-i}, \qquad n\ge 1. \end{equation}
We write  $c=(C_n)_{n\ge 0}$ and then $\Vert c\Vert_{1, {1\over 4^n}}=2$ and $C(z)=Z(c)(z)$ for   $z\in D(0,{1\over 4})$. We may interpret the equality (\ref{cata}) in terms of convolution product in the following closed form
 $$\delta_1* c^{*1}-c+\delta_0=0,$$
 where we deduce that
 \begin{equation}\label{inverse}
 c^{-1}= \delta_0-\delta_1\ast c.
 \end{equation}

\begin{definition}\label{suc} Given the Catalan triangle numbers $(B_{n,k})_{n,k}$ and $(A_{n,k})_{n,k}$ considered in Section 2, we define the Catalan triangle sequences $a_k$ and $b_k$  by
$$
a_k(n):= A_{n+k-1,k}, \qquad b_k(n):= B_{n+k,k}, \qquad n\ge 0,
$$
for $k\ge 1$. Note that $a_1(n)=A_{n,1}=C_n$ and $b_1(n)=B_{n+1,1}=C_{n+1}$ for $n\ge 0$.
\end{definition}

\begin{proposition} \label{cons} For $k\ge 1$, consider the sequences $a_k$ and $b_k$ given in Definition \ref{suc}. Then
\begin{itemize}
\item[(i)] $a_k, b_k\in \ell^1(\N^*, {1\over 4^n})$ and
$$
\Vert a_k\Vert_{1, {1\over 4^n}}=2^{2k-1},\qquad \Vert b_k\Vert_{1, {1\over 4^n}}=2^{2k}.
$$
\item[(ii)] $Z(a_k)(z)=(C(z))^{2k-1}$ and  $Z(b_k)(z)=(C(z))^{2k}$ for $ z\in D(0,{1\over 4})$.
\item[(iii)] $a_k=c^{\ast (2k-2)}$ and $b_k=c^{\ast (2k-1)}$.
\end{itemize}
\end{proposition}
\begin{proof} The item (i) is a consequence of Lemma \ref{tec}. To check (ii), note that
\begin{eqnarray*}
Z(a_k)(z)&=& \sum_{n=0}A_{n+k-1,k}z^n= z^{-k+1}\sum_{m=k-1} A_{m,k}z^{m}=C^{2k-1}(z), \cr
Z(b_k)(z)&=& \sum_{n=0}B_{n+k,k}z^n= z^{-k}\sum_{m=k} B_{m,k}z^{m}=C^{2k}(z), \cr
\end{eqnarray*}
where we have applied fomulae (\ref{gene2}) and (\ref{gene3}). The item (iii) is a straightforward consequence of (ii).
\end{proof}

We may get an alternative proof of Proposition \ref{cons} (i) from item (ii). Note that
\begin{eqnarray*}
\Vert a_k\Vert_{1, {1\over 4^n}}=Z(a_k)({1\over 4})=(C(z))^{2k-1}({1\over 4})=2^{2k-1},\cr
\Vert b_k\Vert_{1, {1\over 4^n}}=Z(b_k)({1\over 4})=(C(z))^{2k}({1\over 4})=2^{2k},
\end{eqnarray*}
 for $k\ge 1$.
\medskip

\begin{proposition} \label{spec} The spectra of the Catalan triangle sequences $(a_k)_{k\ge 1}$ and $(b_k)_{k\ge 1}$ in the algebra $\ell^1(\N^0, {1\over 4^n})$ are given by
$$
\sigma(a_k)=\left(C(\overline{D(0,{1\over 4})})\right)^{2k-1}, \qquad
\sigma(b_k)=\left(C(\overline{D(0,{1\over 4})})\right)^{2k},
$$
for $k\ge 1$. Their boundary is given by
\begin{eqnarray*}
\partial (\sigma(a_k))&=&\left\{2^{2k-1}e^{-i(2k-1)\theta }\left(1-\sqrt{2\vert\sin({\theta\over 2})}\vert e^{i(\pi-\theta)\over 4}\right)^{2k-1}\,\, :\,\,\theta\in (-\pi, \pi)\right\},\\
\partial (\sigma(b_k))&=&\left\{2^{2k}e^{-i2k\theta }\left(1-\sqrt{2\vert\sin({\theta\over 2})}\vert e^{i(\pi-\theta)\over 4}\right)^{2k}\,\, :\,\,\theta\in (-\pi, \pi)\right\}.
\end{eqnarray*}

\end{proposition}

\begin{proof}

As the algebra $\ell^1(\N^0, {1\over 4^n})$  has identity,  the spectrum of an element equals the range of its Gelfand transform (\cite[Theorem 3.4.1]{La}). Moreover as $\sigma(c)=C(\overline{D(0,{1\over 4})})$ (\cite[Proposition 3.2]{[MR3]}), we apply Proposition  \ref{cons} (ii)  to get both first equalities, i.e,
\begin{eqnarray*}
\sigma(a_k)&=&Z(a_k)(\overline{D(0,{1\over 4}}))=\left(C(\overline{D(0,{1\over 4})})\right)^{2k-1}, \\
\sigma(b_k)&=&Z(a_k)(\overline{D(0,{1\over 4})})=\left(C(\overline{D(0,{1\over 4})})\right)^{2k},
\end{eqnarray*}
for $k\ge 1$. As
$$\partial (\sigma(c))=\left\{2e^{-i\theta}\left(1-\sqrt{2\vert\sin({\theta\over 2})}\vert e^{i(\pi-\theta)\over 4}\right)\,\, :\,\,\theta\in (-\pi, \pi)\right\},$$
see \cite[Proposition 3.2]{[MR3]}, we obtain second equalities from previous ones.
\end{proof}

\begin{remark} In the Figure 1, we plot the sets  $\partial (\sigma(c)), \partial (\sigma(b_1))$ and $\partial (\sigma(a_2))$.

\begin{figure}
 \centering
   \includegraphics[scale=0.5]{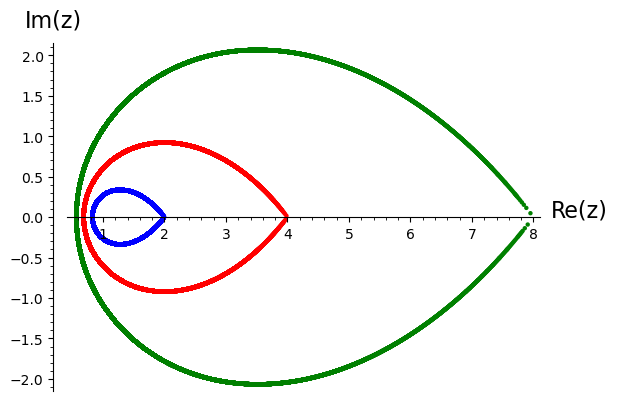}
   \caption{Sets $\partial (\sigma(c))$ in blue, $\partial (\sigma(b_1))$ in red and $\partial (\sigma(a_2))$ in green}
\end{figure}

\end{remark}

Catalan polynomials are defined by the following linear recurrence relation
\begin{equation}\label{recurre22}
{\mathcal P}_{k+2}(z)= {\mathcal P}_{k+1}(z)-z{\mathcal P}_k(z), \qquad k\ge 2,
\end{equation}
and the starting values $\P_0(z)=\P_1(z)=1$. The first values obtained are $\P_2(z)=1-z$, $\P_3(z)=1-2z$ and $\P_4(z)=1-3z+z^2$. The closed form of $\P_k$ is given by the formula
$$
\P_k(z)={(1+\sqrt{1-4z})^{k+1}-(1-\sqrt{1-4z})^{k+1}\over 2^{k+1}\sqrt{1-4z}},
$$
for $k\ge 0$. The bivariate generating function is
$$
{1\over 1-t+zt^2}=\sum_{k\ge 0}\P_k(z)t^k,
$$
see these and other properties in  \cite{Larcombe}. Other interesting property of Catalan polynomials is the following
$$
{d \P_k(z)\over d z}={-1\over 2^{k-1}}\sum_{l=0}^{k-2}(l+2)2^{l}\P_l(z), \qquad k\ge 2,
$$
(\cite[Identity II]{[CLF]}). By induction method, we conclude that the coefficients of $\P_k(z)$ has  alternative signs.

In the next results, we use the usual notation $P(\delta_1)$ where
$$
P(\delta_1):=\sum_{k=0}^n a_k \delta_1^{\ast k}=\sum_{k=0}^n a_k \delta_{ k}
$$
and $P$ is the polynomial, $P(z)=\sum_{k=0}^n a_k z^{ k}$.

\begin{lemma}\label{newwss} Take the Catalan sequence polynomials $(\P_k)_{k\ge 0}$. Then $\P_k(\delta_1)\in \ell^1(\N^0, {1\over 4^n})$,  $\Vert \P_0(\delta_1)\Vert_{1, {1\over 4^n}}=1$ and
$$ \Vert \P_k(\delta_1)\Vert_{1, {1\over 4^n}}=\P_k({-1\over 4})={\alpha_k\over 4^{k-1}}, \qquad k\ge 1,
$$
where $\alpha_1=1$, $\alpha_2=5$ and $\alpha_k=4(\alpha_{k-1}+\alpha_{k-2})$ for $k\ge 3$.
\end{lemma}
\begin{proof} It is clear that  $\P_k(\delta_1)\in \ell^1(\N^0, {1\over 4^n})$ and $\Vert \P_0(\delta_1)\Vert_{1, {1\over 4^n}}=1$. As  the coefficients of polynomials $(\P_k)_{k\ge 0}$ have alternative  signs, we conclude that
\begin{eqnarray*}
\Vert \P_k(\delta_1)\Vert_{1, {1\over 4^n}}&=&\sum_{j=0}^k a_j \left({ -1 \over  4}\right)^{ j}=\P_k({-1\over 4})\cr&=&{(1+\sqrt{2})^{k+1}-(1-\sqrt{2})^{k+1}\over \sqrt{2}2^{k+1}}={\alpha_k\over 4^{k-1}},\cr
\end{eqnarray*}
where the integer sequence $(\alpha_k)_{k\ge 1}$ is numbered as A086347 in \cite{[Sl]} and treated in detail there.
\end{proof}

\begin{remark}{\rm The first values of the sequence $(\alpha_k)_{k\ge 1}$ are $1,$ $5,$ $24,$ $116,$ $560...$. This sequence is an example of
generalized Fibonacci numbers $g(k) = cg(k-1) + dg(k-2)$ for $k\ge 2$ and  seed values $g(0)=a$ and $g(1)=b$ ($a,b,c,d\in \N$.) }

\end{remark}
\begin{theorem}\label{convoinve} For $k \ge 1$,
$$
(c^{\ast k})^{-1}= \P_{k+1}(\delta_1)+(-c\ast\delta_1)\ast \P_{k}(\delta_1).
$$
Moreover $\Vert (c\ast c)^{-1}\Vert_{1, {1\over 4^n}} ={3\over 2}$ and  $\Vert (c^{\ast k})^{-1}\Vert_{1, {1\over 4^n}}\le{1\over 4^{k}}\left(\alpha_{k+1}+2\alpha_{k}\right)$ for $k\ge 1$, where $(\alpha_k)_{k\ge 1}$ are defined in Lemma \ref{newwss}.
\end{theorem}

\begin{proof} Note that $c^{-1}= \delta_0-\delta_1\ast c$, see formula (\ref{inverse}) and then
\begin{eqnarray*}
(c\ast c)^{-1}&=& c^{-1}\ast c^{-1}= \delta_0- 2\delta_1\ast c +\delta_1\ast(\delta_1\ast c\ast c)\cr
&=& \delta_0-\delta_1-\delta_1\ast c= \P_2(\delta_1)+(-c\ast \delta_1)\ast \P_1(\delta_1),
\end{eqnarray*}
where we have applied that $\delta_1\ast c^{\ast 1}= c-\delta_0$. By induction, we have that
\begin{eqnarray*}
(c^{\ast(k+1)})^{-1}&=& c^{-1}\ast (c^{\ast k})^{-1}= (\delta_0-\delta_1 \ast c)\ast (\P_{k+1}(\delta_1)+(-c\ast\delta_1)\ast \P_{k}(\delta_1))\cr&=&
\P_{k+1}(\delta_1)-\delta_1\ast c\ast \P_{k+1}(\delta_1)- \delta_1\ast \P_{k}(\delta_1)\cr&=& \P_{k+2}(\delta_1)+ (-c\ast \delta_1)\ast \P_{k+1}(\delta_1),
\end{eqnarray*}
where we have applied the recurrence relation (\ref{recurre22}).

Finally, we apply Lemma  \ref{newwss} to get
$$
\Vert (c^{\ast k})^{-1}\Vert_{1, {1\over 4^n}}\le\Vert \P_{k+1}(\delta_1)\Vert_{1, {1\over 4^n}}+{1\over 2}\Vert \P_{k}(\delta_1)\Vert_{1, {1\over 4^n}}=  {1\over 4^{k}}\left(\alpha_{k+1}+2\alpha_{k}\right)
$$
for $ k\ge 1$. \end{proof}

\section{ Powers of Catalan generating functions for bounded operators}

In this section, we consider the particular case that $T$ is a linear and bounded operator on the Banach space $X$, $T\in {\mathcal B}(X)$,   such that
\begin{equation}\label{pb}
\sup_{n\ge 0} {\Vert 4^n T^n\Vert}:=M<\infty,
\end{equation}
i.e., ${\displaystyle{4T}}$ is a power-bounded operator. In this case $\sigma(T)\subset \overline{D(0,{1\over 4})}$. Under the condition (\ref{pb}), we define the Catalan generating function, $C(T)$, by
\begin{equation} \label{serie}
C(T):=\sum_{n\ge 0}C_nT^n,
\end{equation}
see \cite[Section 5]{[MR3]}. The bounded operator $C(T)$ may be considered as the  image of the Catalan sequence $c=(C_n)_{n\ge 0}$ in the algebra homomorphism $\Phi: \ell^1(\N^0, {1\over 4^n}) \to {\mathcal B}(X)$ where
$$
\Phi(a)x:=\sum_{n\ge 0} a_nT^n(x), \qquad  a=(a_n)_{n\ge 0}\in  \ell^1(\N^0, {1\over 4^n}) , \quad x\in X,
$$
i.e., $\Phi(c)=C(T)$. The $\Phi$ algebra homomorphism (also called functional calculus) is presented in some functional analysis textbooks, for example \cite[Chapter 13 and 14]{Mu}.

\begin{theorem}\label{genne} Given  $T\in {\mathcal B}(X)$ such that ${4T}$ is power-bounded and $c=(C_n)_{n\ge 0}$ the Catalan sequence. Then
\begin{itemize}
\item[(i)] The powers $(C(T))^{2k-1}= \Phi(a_k)$ and $(C(T))^{2k}=\Phi(b_k)$ for $k\ge 1$,
and
$$
\Vert (C(T)^j\Vert \le (C(\Vert T\Vert))^j, \qquad j\ge 1.
$$
\item[(ii)] The operator $C(T)$ is invertible, $(C(T))^{-1}= I-TC(T)$,
$$
(C(T))^{-(j+1)}= \P_j(T)-TC(T) \P_{j-1}(T)\qquad j\ge 1,
$$
 $\Vert C(T)^{-1}\Vert\le 1+{1\over 2}\sup_{n\ge 0} {\Vert 4^n T^n\Vert}$, $\Vert C(T)^{-2}\Vert\le {3\over 2}\sup_{n\ge 0} {\Vert 4^n T^n\Vert}$ and
$$
\Vert(C(T))^{-(j+1)}\Vert\le {1\over 4^{j}}\sup_{n\ge 0} {\Vert 4^n T^n\Vert}\left(\alpha_{j+1}+2\alpha_{j}\right), \qquad j\ge 1, $$
where $(\alpha_j)_{j\ge 1}$ are defined in Lemma \ref{newwss}.

\item[(iii)] Take $(P_n)_{n\ge 0}$ and $(Q_n)_{n\ge 0}$ polynomials given in Definition \ref{polys}. Then
\begin{eqnarray*}
\sum_{n\ge 0}P_n(z)T^n&=& {C(T)-(z+1)I\over T(1+z)^2-zI},\cr
\sum_{n\ge 0}Q_n(z)T^n&=& {(C(T)-(z+1)I)(z+1)\over T(1+z)^2-zI},
\end{eqnarray*}
for $\vert z\vert <1$.
\item[(iv)] The spectral mapping theorem holds for $(C(T))^n$, i.e, $\sigma((C(T))^n)=C^n(\sigma(T))$ for $n\in \Z.$
\end{itemize}

\end{theorem}
\begin{proof} (i) From (\ref{serie}),   $\Phi(c)=C(T)\in {\mathcal B}(X)$ as we have commented above. By Proposition  \ref{cons} (iii), we have
\begin{eqnarray*}
(C(T))^{2k-1}&=&(\Phi(c))^{2k-1}=\Phi(c^{\ast(2k-2)})=\Phi(a_k),\cr
(C(T))^{2k}&=&(\Phi(c))^{2k}=\Phi(c^{\ast(2k-1)})=\Phi(b_k),
\end{eqnarray*}
for $k\ge 1$. By Proposition  \ref{cons} (ii), we get
\begin{eqnarray*}
\Vert (C(T))^{2k-1}\Vert = \Vert\Phi(a_k)\Vert\le \sum_{j\ge 0}a_{k}(j)\Vert T\Vert^j=(C(\Vert T\Vert))^{2k-1},\cr
\Vert (C(T))^{2k}\Vert = \Vert\Phi(b_k)\Vert\le \sum_{j\ge 0}b_{k}(j)\Vert T\Vert^j=(C(\Vert T\Vert))^{2k},
\end{eqnarray*}
for $k\ge 1$ and we conclude the proof of (i).

\noindent (ii) As the homomorphism $\Phi$ is continuous, we apply the formula (\ref{inverse}) to get
$$
C(T)(I-TC(T))=\Phi(c)(\Phi(\delta_0-\delta_1\ast c))=\Phi(c-\delta_1\ast c^{\ast 1})=\Phi(\delta_0)=I.
$$
In fact $(C(T))^{-1}= \Phi(c^{-1})$ and
$$
(C(T))^{-(j+1)}= \Phi((c^{-1})^{\ast j})=\Phi(c^{\ast j})^{-1}= \P_{j+1}(T)-TC(T) \P_{j}(T),\qquad j\ge 1,
$$
where we have applied Theorem \ref{convoinve} and $\Phi$ is an algebra homomorphism. The estimation of $\Vert (C(T))^{-(j+1)}\Vert $ follows also from Theorem \ref{convoinve}.

\noindent(iii) We follow similar ideas to those shown in Theorem \ref{generating} and we check
$$
\sum_{n\ge 0}P_n(z)T^n= {C(T)-(z+1)I\over T(1+z)^2-zI}, \qquad \displaystyle{\sum_{n\ge 0}Q_n(z)T^n= {(C(T)-(z+1)I)(z+1)\over T(1+z)^2-zI}},
$$
 for $\vert z\vert <1.$

\noindent(iv) Since ${4T}$ is power bounded, the spectral mapping theorem for $C^n(T)$ may found in \cite[Theorem 2.1]{[Dun]} and then
$\sigma((C(T))^n)=C^n(\sigma(T))$ for $n\in \Z.$
\end{proof}

\begin{remark} As $\sigma(T)\subset \overline{D(0,{1\over 4})}$, we apply Proposition \ref{spec} to conclude that
$$\sigma(C^n(T))\subset C^n(\overline{D(0,{1\over 4})}), \qquad n\in \Z.
$$
\end{remark}
%%%%%%%%%%%%%%%%%%%%%
%%%%%%%%%%%%%%%%%%%%%
\section{Examples, applications and final comments}

 In this section we present some particular examples of operators $T$  for which we solve  the equation  (\ref{Ceq}), calculate $C(T)$ and $(C(T))^k$ for $k\in \Z$. In the subsection 5.1, we consider the Euclidean space $\C^2$ and some matrices $T$. To  resolve this matrix equation, we need to solve a system of four quadratic equations.  We also calculate $(C(T))^{n}$ for these matrices. In  subsection 5.2  we check $C(a)$ for some $a\in \ell^1(\N^0, {1\over 4^n})$. Finally we present some ideas to continue this research in subsection 5.3.

\subsection{Matrices on $\C^2$}\label{matrix} We consider the Euclidean space $\C^2$ and the operator $T=\begin{pmatrix}
\lambda & 0\\
0 & \mu
\end{pmatrix},
$
 with $0\not=\lambda, \mu\in \C$. For $\lambda=\mu$, the solution is presented in \cite[Subsection 6.1]{[MR3]}. For $\lambda\not=\mu$,  the solution of  (\ref{Ceq}) is given by
$$
Y=\begin{pmatrix}
{1\pm\sqrt{1- 4\lambda }\over 2\lambda} & 0\\
 0 & {1\pm\sqrt{1-4\mu}\over 2\mu}
\end{pmatrix}, \qquad
$$
 where  the allowed signs are all four combinations. In the case that $\vert \lambda \vert, \vert \mu \vert\le {1\over 4}$, note that
$$
(C(T))^j=\begin{pmatrix}
(C(\lambda))^j & 0\\
0 & (C(\mu))^j
\end{pmatrix}. \qquad
$$
for $j\in \Z$.

\medskip
Now we study the case $T=\begin{pmatrix}
0 & \lambda\\
\lambda & 0
\end{pmatrix}$ with $\lambda \in \C\backslash\{0\}$. When $\vert \lambda \vert\le {1\over 4}$, we get that
\begin{equation}\label{matt}
C(T)=\begin{pmatrix}
C_e(\lambda) & C_o(\lambda)\\
 C_o(\lambda) &  C_e(\lambda)
\end{pmatrix},
\end{equation} where functions $C_e $ and $C_o$ are functions given by
\begin{eqnarray*}
C_e(\lambda):=\sum_{n=0}^\infty C_{2n}\lambda^{2n}&=&{\sqrt{1+4\lambda}-\sqrt{1-4\lambda}\over 4\lambda},\cr
C_o(\lambda):=\sum_{n=0}^\infty C_{2n+1}\lambda^{2n+1}&=&{2-\sqrt{1+4\lambda}-\sqrt{1-4\lambda}\over 4\lambda}.\cr
\end{eqnarray*}
Note that $C(T)$ is one of the four solutions of (\ref{Ceq}), see \cite[Section 6.1]{[MR3]}.  As
\begin{eqnarray*}
\begin{pmatrix}
a & b\\
 b &  a
\end{pmatrix}^{2n}&=&\left(a^{2n}+{2n\choose 2}a^{2(n-2)}b^2+\dots...+b^{2n} \right)\begin{pmatrix}
1 & 0\\
0& 1
\end{pmatrix}\cr&+&\left({2n\choose 1}a^{2n-1}b+\dots...+{2n\choose 1}ab^{2n-1} \right)\begin{pmatrix}
0 & 1\\
1 & 0
\end{pmatrix},\cr
\begin{pmatrix}
a & b\\
 b &  a
\end{pmatrix}^{2n+1}&=&\left(a^{2n+1}+\dots...+{2n+1\choose 2n}ab^{2n} \right)\begin{pmatrix}
1 & 0\\
0& 1
\end{pmatrix}\cr&+&\left({2n+1\choose 1}a^{2n}b+\dots...+b^{2n+1} \right)\begin{pmatrix}
0 & 1\\
1 & 0
\end{pmatrix},\cr
\end{eqnarray*}
we use (\ref{matt}) to get new generating formulae for Catalan triangle numbers.

\begin{theorem}Take $n\ge 0$ and $z\in \overline{D(0,{1\over 4})}$. Then

\begin{eqnarray*}
\sum_{k=n}^\infty B_{2k-n,n }z^{2k}&=&z^{2n}\left( C_e^{2n}(z)+{2n\choose 2}C_e^{2(n-2)}(z)C_o^2(z)+\dots...+C_o^2(z)\right),\cr
\sum_{k=n}^\infty B_{2k+1-n,n }z^{2k+1}&=&z^{2n}\left({2n\choose 1} C_e^{2n-1}(z)C_o(z)+\dots...+{2n\choose 1}C_e(z)C_o^{2n-1}(z) \right),\cr
\sum_{k=n}^\infty A_{2k-1-n,n }z^{2k}&=&z^{2n}\left(C_e^{2n-1}(z)+\dots...+{2n-1\choose 2n-2}C_e(z)C_o^{2n-2}(z) \right),\cr
\sum_{k=n}^\infty A_{2k-n,n }z^{2k+1}&=&z^{2n}\left({2n-1\choose 1}C_e^{2n-1}(z)C_o(z)+\dots...+C_o^{2n-1}(z) \right),\cr
\end{eqnarray*}
\end{theorem}
\begin{proof} Take $\vert z \vert\le {1\over 4}$ and we consider $T=\begin{pmatrix}
0 & z\\
z & 0
\end{pmatrix}$. We apply  Theorem \ref{genne} to get
\begin{eqnarray*}
C(T)^{2n}&=&\sum_{j=0}^\infty b_n(j)T^j=\sum_{l=0}^\infty b_n(2l)z^{2n}I +\sum_{l=0}^\infty b_n(2l+1)z^{2n+1}\begin{pmatrix}
0 & 1\\
1 & 0
\end{pmatrix} \cr
&=&\sum_{l=0}^\infty B_{2l+n,n}z^{2n}I +\sum_{l=0}^\infty B_{2l+1+n,n}z^{2n+1}\begin{pmatrix}
0 & 1\\
1 & 0
\end{pmatrix} \cr
&=&{1\over z^{2n}}\left(\sum_{k=n}^\infty B_{2k-n,n}z^{2k}I +\sum_{k=n}^\infty B_{2k+1-n,n}z^{2k+1}\begin{pmatrix}
0 & 1\\
1 & 0
\end{pmatrix}\right),
\end{eqnarray*}
and we conclude the first two equalities. Similarly, we consider $C(T)^{2n+1}$ and show the second two equalities.
\end{proof}
\medskip
Finally we study the case $T=\begin{pmatrix}
\lambda & \mu\\
0 & \lambda
\end{pmatrix}$ with $\lambda, \mu \in \C\backslash\{0\}$. The solutions of (\ref{Ceq}) are given by
$$
Y=\begin{pmatrix}
{a} & {\mu(a-1)\over \lambda (1-2\lambda a)}\\
0 & a
\end{pmatrix}
$$
where $a$ is any solution of the quadratic Catalan equation $\lambda a^2-a+1=0$. In the case that $\vert \lambda \vert\le {1\over 4}$, we get that
$$
C(T)=\begin{pmatrix}
C(\lambda)& {\mu((C(\lambda)-1)\over \lambda (1-2\lambda C(\lambda))}\\
 0 &  C(\lambda)
\end{pmatrix},
$$
and
$$
(C(T))^j=\begin{pmatrix}
(C(\lambda))^j& n(C(\lambda))^{j-1}{\mu((C(\lambda)-1)\over \lambda (1-2\lambda C(\lambda))}\\
 0 &  (C(\lambda))^j
\end{pmatrix},
$$
for $j\ge 1$. As $
(C(T))^{-1}={1\over (C(\lambda))^2}\begin{pmatrix}
C(\lambda)& -{\mu((C(\lambda)-1)\over \lambda (1-2\lambda C(\lambda))}\\
 0 &  C(\lambda)
\end{pmatrix}
$, we get that
$$
(C(T))^{-j}={1\over (C(\lambda))^{2j}}\begin{pmatrix}
(C(\lambda))^j& -n(C(\lambda))^{j-1}{\mu((C(\lambda)-1)\over \lambda (1-2\lambda C(\lambda))}\\
 0 &  (C(\lambda))^j
\end{pmatrix},
$$
for $j\ge 1$.

\subsection{Catalan operators on  $\ell^p$} We consider the space of sequences $\ell^p(\N^0, {1\over 4^n})$  where
$$
\Vert a\Vert_{p, {1\over 4^n}}:=\left(\sum_{n=0}^\infty {\vert a_n\vert^p\over 4^{np}}\right)^{1\over p}<\infty,
$$
for $1\le p<\infty$ and  $\ell^\infty(\N^0, {1\over 4^n})$ the space of sequences embedded with the norm
$$
\Vert a\Vert_{\infty, {1\over 4^n}}:=\sup_{n\ge 0}{\vert a_n\vert\over 4^{n}}<\infty.
$$
 Note that $\ell^1(\N^0, {1\over 4^n})\hookrightarrow \ell^p(\N^0, {1\over 4^n})\hookrightarrow \ell^\infty(\N^0, {1\over 4^n})$.

 Now we consider sequences $c, (a_k), (b_k)\in \ell^1(\N^0, {1\over 4^n})$,  the Catalan triangle sequences given in Definition \ref{suc}   and convolution operators $C(f):= c\ast f$, $C^{2k}(f)= b_k\ast f$ and $C^{2k-1}(f)= a_k\ast f$  for $f\in \ell^p(\N^0, {1\over 4^n})$ and $k\ge 1$ with $1\le p\le \infty$. By Theorem \ref{genne} (iv), we get that
 $$
 \sigma(C^n)=C^{n}(\sigma(\delta_1))=C^n(\overline{D(0,{1\over 4})}),\qquad n \ge 1.
 $$
 Note that  the set $\sigma(C^n)$ independent on $p$ and coincides with the spectrum of the power of Catalan sequence $c$ in  $\ell^1(\N^0, {1\over 4^n})$
 (Proposition \ref{spec}).

\subsection{A future research}

Given $a, b\not=0\ \in \C$, the quadratic equation
\begin{equation}\label{genequat}
{bz\over 2}y^2-y+{a\over 2b}=0,
\end{equation}
has two solutions given by
$$
y={1\pm \sqrt{1-za}\over bu}.
$$
We define $\displaystyle{C^{a,b}(z):= {1- \sqrt{1-za}\over bz}};$
note that
$
\displaystyle{C^{a,b}(z)={a\over 2b}C({az\over 4})}
$ and
$$
C^{a,b}(z)=\sum_{n\ge 0}{a^{n+1}\over 2^{2n+1}b}C_n z^n.
$$

It would be natural to consider a vector-valued version of equation  (\ref{genequat}) for $a, b, z\in {\mathcal B}(X)$.

%%%%%%%%%%%%%%

\noindent{\bf Acknowledgment.} We thank the referee for his/her very careful review of this paper, and for the comments, corrections and
suggestions that ensued. A major version of the paper  has been carried out to take them
into account and the paper has been significantly improved.

\noindent. This work does not have any conflicts of interest.

%%%%%%%%%%%%%%
%%%%%%%%%%%%%%%

\end{document}